\documentclass[authoryear,preprint]{elsarticle}

\usepackage{graphicx,amssymb,amsfonts,amsmath,latexsym,mathtools}
\usepackage[round]{natbib}

\journal{ ...}

\newtheorem{theo}{Theorem}

\newtheorem{coro}{Corollary}

\newcommand{\epr}{\hfill $\Box$\mbox{}\\ }
\newcommand{\bpr}{\noindent {\bf Proof.} \hspace{1 em}}

\newcommand\floor[1]{\lfloor#1\rfloor}

\newcommand\bcar[1]{\bigg\lvert#1\bigg\rvert}

\newcommand{\RR}{\mathbb{R}}
\newcommand{\NN}{\mathbb{N}}
\newcommand{\ZZ}{\mathbb{Z}}

\newcommand{\m}{\overline{m}}

\newcommand{\bigO}{\mathcal{O}}

\newcommand{\coef}[1]{\ensuremath{c_{#1}}}

\DeclarePairedDelimiter\abs{\lvert}{\rvert} 
\DeclarePairedDelimiter\card{\lvert}{\rvert}
\makeatletter
\let\oldabs\abs
\def\abs{\@ifstar{\oldabs}{\oldabs*}}
\let\oldcard\card
\def\card{\@ifstar{\oldcard}{\oldcard*}}
\makeatother

\begin{document}

\begin{frontmatter}

\title{The size of the largest antichains in  products of linear orders}
\tnotetext[t1]{Authors are listed alphabetically and have equally contributed.}

\author[Lamsade]{D. Bouyssou} \ead{bouyssou@lamsade.dauphine.fr}
\author[Ghent]{T. Marchant\corref{cor1}} \ead{thierry.marchant@ugent.be}
\cortext[cor1]{Corresponding author} 
\author[Mons]{M. Pirlot} \ead{marc.pirlot@umons.ac.be}

\address[Lamsade]{LAMSADE, UMR 7243, CNRS, Universit\'e Paris-Dauphine, PSL Research University, 75016 Paris, France}
\address[Ghent]{Ghent University, Dunantlaan 1, 9000 Ghent, Belgium} 
\address[Mons]{Universit\'e de Mons, rue de Houdain 9, 7000 Mons, Belgium} 

\begin{abstract}
   We present new exact and asymptotic results about the size of the largest antichain in the product of $n$ linear orders.
\end{abstract}


\end{frontmatter}

\section{Introduction}

Antichains in the poset $\{0,1\}^{n}$ equipped with the standard partial ordering are well-studied and have many different interpretations \citep[e.g.][]{ErsekUyanikSMP17}.
An expression for the maximal size of such antichains in $\{0,1\}^{n}$  is given by a classical  theorem of  \cite{Sperner1928}. 
If we consider the more general poset $\{1, \ldots, m\}^{n}$ also equipped with the standard partial ordering, an expression for the size of the largest antichain is given in \cite{Sander1993}. Sander also provides  asymptotic results when $m$ is fixed and $n$ goes to infinity. The interest of Sander in this problem arose form a recreational mathematics problem posed in \cite{Motek1986}. Actually, antichains  and, hence, maximal antichains, in the poset $\{1, \ldots, m\}^{n}$ are of interest in many domains. For instance in game theory, \cite{HsiaoRaghavan1993} define a multichoice cooperative game as a real-valued mapping on $\{1, \ldots, m\}^{n}$, where $n$ is the number of players and $\{1, \ldots, m\}$ denotes the set of ordered actions that each player can take. A profile in such a game is  a vector $x=(x_{1}, \ldots, x_{n}) \in \{1, \ldots, m\}^{n}$ and represents the actions taken by each agent. A winning profile is such that the value of the game at that profile is 1. A winning profile $x$ is minimal if there is no other winning profile $y$ such that $y \leq x$. If a game is monotone, then the set of all minimal winning profiles is an antichain.
Besides, \cite{Grabisch2016} shows that antichains in $\{1, \ldots, m\}^{n}$ play an important role in the analysis of these multichoice games. 
Our personal interest in antichains in the poset $\{1, \ldots, m\}^{n}$ stems from the analysis of a multicriteria sorting model presented in \cite{FernandezFigueiraNavarroRoy2017}. In this context, the size of maximal antichains corresponds to the maximum number of profiles needed to represent a twofold ordered partition in their model, whenever such a representation is possible. 
Another paper about antichains in $\{1, \ldots, m\}^{n}$ is \cite{Tsai2018}: it presents an upper bound for the number of antichains (a generalization of Dedekind numbers).

In the present paper, we will extend Sander's results in two directions. First, we will present an exact expression for the size of the largest antichain in the heterogeneous product $\Pi_{i=1}^{n} \{1, \ldots, m_{i}\}$. Then, we will provide asymptotic results for the size of the largest antichain in $\{1, \ldots, m\}^{n}$ when $n$ is fixed and $m$ goes to infinity.

\section{Notation and definitions}
\label{sec:notation}

Let $P$ be a set and $\leq$ be a binary relation defined on $P$, satisfying (i) reflexivity ($\forall x \in P, x \leq x$), (ii) antisymmetry ($\forall x,y \in P, x \leq y$ and $y \leq x \iff x=y$) and (iii) transitivity ($\forall x,y,z \in P, x \leq y$ and $y \leq z \Rightarrow x \leq z$). The pair $(P, \leq)$ is called a partially ordered set (\emph{poset})\footnote{Most definitions about posets are taken from \cite{ProctorSaksSturtevant1980}.}. When there is no ambiguity, the poset $(P, \leq)$ is simply denoted by $P$. 
For all $x,y$ belonging to a poset $P$, we say that $x$ and $y$ are comparable if $x \leq y$ or $y \leq x$. A chain of $P$ is a totally ordered subset of $P$. A \emph{linear order} on $P$ is a poset such that $P$ is a chain.
An \emph{antichain} of $P$ is a subset of pairwise incomparable elements.  A largest antichain is an antichain of maximal cardinality.

Let $(P, \leq_{P})$ and $(Q, \leq_{Q})$ be two posets. The product poset $(P \times Q, \leq)$ is defined to be the set of all pairs $(a,b), a \in P, b \in Q$, with the order given by $(a, b) \leq (a', b')$ if and only if $(a \leq_{P} a')$ and $(b \leq_{Q} b')$. Let  $n$ be a positive integer and  $\overline{m} = (m_{1}, \ldots, m_{n})$ be an element of $\NN^{n}$, where $\NN$ denotes the set of positive integers. For any $a \in \NN$, let $[a]$ denote the set $\{1, \ldots, a\}$. For any $i \in [n]$, the  poset $([m_{i}], \leq)$,  where $\leq$ is the usual ordering of the natural numbers, is a linear order (also called a chain). The product of these $n$ linear orders is the poset $(\Pi_{i=1}^{n} [m_{i}], \preceq)$ where $\preceq$ is defined as follows: for all $x,y \in \Pi_{i=1}^{n} [m_{i}]$, $x \preceq y$ iff $x_{i} \leq y_{i}$ for all $i \in [n]$. The poset $(\Pi_{i=1}^{n} [m_{i}], \preceq)$ is also called a direct product of chains \citep{CaspardLeclercMonjardet2012}. When $\m$ is such that $m_{i}=m$ for all $i \in [n]$, then the Cartesian product $\Pi_{i=1}^{n} [m_{i}]$ is homogeneous and can be written as $[m]^{n}$.

The size of the largest antichains in $\Pi_{i=1}^{n} [m_{i}]$ and $[m]^{n}$ will  respectively be denoted by $s(\m)$ and $S(m,n)$.
\cite{Sperner1928} has proved that the size of the largest antichain in $[2]^{n}$ is 
$$S(2,n) = \binom{n}{\floor{n/2}}.$$
When $n$ is large, a convenient approximation for $S(2,n)$  is obtained using Stirling's formula:
$S(2,n) \sim 2^{n} \sqrt{2/\pi n}.$
Later, \cite{Sander1993} has proved that  the size of the largest antichain in $[m]^{n}$ is 
\begin{equation}
\label{eq:Sander}
S(m,n) = \sum_{j=0}^{ \lfloor g/m \rfloor } (-1)^{j} \binom{n}{j} \binom{n-1+g-mj}{n-1},
\end{equation}
with $g=\lfloor n(m-1)/2 \rfloor$. Sander has also provided a bound\footnote{This bound is later rediscovered by \cite{MattnerRoos2008} in a different context. In addition, \cite{MattnerRoos2008} note that \eqref{eq:Sander} can be found in \cite{Moivre1756} as the solution of a probability problem.} and some asymptotic results for $S(m,n)$ when $m$ is fixed. Notice that $S(m,n)$ corresponds to Sequence A077042 in the Online Encyclopedia of Integer Sequences \citep{OEISB}.

Section~\ref{sec.hetero} is devoted to the general case of heterogeneous products and   presents some exact results about $s(\m)$. In Section~\ref{sec.homo}, we  consider the special case of homogeneous products and we  present a new exact result about $S(m,n)$ and also an asymptotic result when $n$ is fixed. 

\section{Heterogeneous product}
\label{sec.hetero}

Let us define $m_{I}=  \sum_{i \in I} m_{i}$ and 
\begin{equation}
\label{eq.h}
h= \left\lfloor\frac{n+\sum_{i \in N} m_{i}}{2} \right\rfloor. 
\end{equation}
Our  result about heterogeneous products is the following.

\begin{theo}
\label{theo:SanderGeneralized}
For all $\m = (m_{1}, \ldots, m_{n}) \in \NN^{n}$, 
\begin{equation}
\label{eq.Sperner.G}
s(\m) =\sum_{I \subseteq [n]:m_{I}<h-n}  \binom{h-m_{I}-1}{n-1}   (-1)^{\card{I}} .
\end{equation}
\end{theo}
Before proving this result, we recall some definitions and results about posets.
Let $(P, \leq)$ be a poset. For any $x, y \in P$, we say that $y$ \emph{covers} $x$ in $P$ iff $x < y$ and there is no $z$ such that $x<z<y$.
A \emph{ranking} (or grading) of a poset $P$ is a partition of $P$ into (possibly empty) sets $P_{i}$ $(i \in \ZZ)$ such that, for each $i$, every element in $P_{i}$ is covered only by elements   in  $P_{i+1}$. The set $P_{i}$ is called the $i$th rank of $P$. If a poset admits a ranking, then we say that it is ranked (or graded).

The \emph{Whitney numbers} of a ranked poset $P$ are $\{p_{i}: i \in \ZZ \}$, where $p_{i}$ is the cardinality of $P_{i}$. Let $p_{k}$ be a largest Whitney number of a ranked poset; we say $P$ is \emph{rank-unimodal} if $p_{i} \geq p_{i-1}$, for $i  \leq k$ and $p_{i} \geq p_{i+1}$, for $i \geq  k$. The Whitney numbers are said to be symmetric, or $P$ is said to be \emph{rank-symmetric}, if there exists a $d$ such that $p_{i} = p_{d-i}$ for all $i$. Let $P$ and $Q$ be ranked posets and
 $R = P \times Q$. The rankings of $P$ and $Q$ induce the following ranking on $R$:
$$R_{l} = \bigcup_{i} \ (P_{i} \times Q_{l-i}),$$
and $R$ has Whitney numbers 
\begin{equation}
\label{eq.whitney}
r_{l} = \sum_{i  \in \ZZ} p_{i}q_{l-i}.
\end{equation}

 A $k$-family is a subset of $P$ containing no chain of size $k + 1$. 
Equivalently, a $k$-family is a union of $k$ (possibly empty) antichains, so that if $P$ is ranked, any union of $k$ ranks is a $k$-family. 
$P$ is said to be \emph{Sperner} if the rank of largest size is an antichain of maximum size, and $P$ is $k$-Sperner if the union of the $k$ largest ranks is a $k$-family of maximum size. 
$P$ is \emph{strongly Sperner} if it is $k$-Sperner for all $k \geq  1$. 
$P$ is a \emph{Peck} poset if it is strongly Sperner, rank-unimodal, and rank-symmetric. Theorem 3.2 in \cite{ProctorSaksSturtevant1980} shows that  the product of two Peck posets is a Peck poset.

\medskip

\noindent {\bf Proof of Theorem~\ref{theo:SanderGeneralized}.} \hspace{0.5 em}
In this proof, for the sake of brevity, we use $X$ to denote the poset $(\Pi_{i=1}^{n} [m_{i}], \preceq)$ or the set  $\Pi_{i=1}^{n} [m_{i}]$.
For each $i \in [n]$, the poset $([m_{i}], \leq)$ is a Peck poset. Hence, $n-1$ applications of Theorem 3.2 in \cite{ProctorSaksSturtevant1980} show that $X$ is also a Peck poset. Let us consider a ranking of $X$ such that the minimal element $(1, \ldots,1)$ in $X$ has rank $n$, which is the sum of the coordinates of the minimal element.  The maximal element $(m_{1}, \ldots, m_{n})$ has rank $\sum_{i \in [n]} m_{i}$. 
Let $p^{i}_{j}$ be the $j$th Whitney number of $([m_{i}], \leq)$; it is equal to 1 for each non-empty rank.
Since $X$ is rank-unimodal, and rank-symmetric, a maximal Whitney number corresponds to the median rank $h$ (defined by \eqref{eq.h}).
Because of the Sperner property, this rank is also an antichain of maximum size and $n-1$ applications of \eqref{eq.whitney} show that the size of this  antichain is
$$
s(\m) = \sum_{i_{n-1} \in \ZZ} \ldots \sum_{i_{2} \in \ZZ} \sum_{i_{1} \in \ZZ} p^{1}_{i_{1}} p^{2}_{i_{2}}  \ldots p^{n-1}_{i_{n-1}} p^{n}_{h-i_{1}-i_{2}- \ldots - i_{n-1}}.
$$
Define $i_{n}= h-\sum_{j \in [n-1]} i_{j}$
and we obtain
\begin{equation}
\label{eq.whitney.sum}
s(\m) = \sum_{i_{n-1} \in \ZZ} \ldots \sum_{i_{2} \in \ZZ} \sum_{i_{1} \in \ZZ} p^{1}_{i_{1}} p^{2}_{i_{2}}  \ldots p^{n-1}_{i_{n-1}} p^{n}_{i_{n}}
=\sum_{ \substack{ i_{1}+ \ldots + i_{n}=h \\ 1 \leq i_{j} \leq m_{j}, \forall j \in [n] }} 1.
\end{equation}
Hence, an antichain of maximum size in $X$
is the set 
$$A=\{ x \in X : \sum_{i \in [n]} x_{i} = h \}$$ 
and a generating function for $s(\m)$ is defined as follows, for all $x \in \RR$, $\vert x \vert<1$,
$$
f(\m,x) = (x^{1}+ \ldots + x^{m_{1}}) \times \ldots \times (x^{1}+ \ldots + x^{m_{n}}).
$$
We also have 
$$
f(\m,x) = 
\left( \frac{x}{1-x} \right)^{n} \left(  \sum_{I \subseteq [n]} (-1)^{\card{I} } x^{ m_{I}} \right).
$$
 We denote by $\coef{p}(A(x))$ the coefficient of $x^p$ in the polynomial $A(x)$. Hence $s(\m) = \coef{h}(f(\m,x))$. A property of  products of polynomials, which extends to absolutely convergent series, is:
$$\coef{p}(A(x)B(x)) = \sum_{k=0}^{p} \coef{k}(A(x)) \times \coef{p-k}(B(x)). $$
In our case, we have 
\begin{align*}
s(\m) & =  \coef{h}(f(m_{1}, \ldots, m_{n},x)) \\
	 &= \coef{h}\left(\left( \frac{x}{1-x} \right)^{n} \times   \sum_{I \subseteq [n]} (-1)^{\card{I}} x^{ m_{I}} \right)  \\
	& = \sum_{j=0}^{h}  \coef{j}\left( \frac{x}{1-x} \right)^{n} \times  \coef{h-j}\left(\sum_{I \subseteq [n]} (-1)^{\card{I}} x^{ m_{I}}\right)    .
\end{align*}
For $|x|<1$, we have
$$\coef{j}\left( \frac{x}{1-x} \right)^{n} = \coef{j}(x+x^{2}+x^{3}+ \ldots)^{n}$$
and $\coef{j}\left(\frac{x}{1-x} \right)^{n}$ is thus the number of $n$-tuples of positive integers whose sum is equal to $j$. It is $\binom{j-1}{n-1}$ \cite[p.835]{MacMahon1893}.
We also have 
$$
\coef{h-j}\left(\sum_{I \subseteq [n]} (-1)^{\card{I}} x^{ m_{I}}\right)= \sum_{I \subseteq [n]:m_{I}=h-j}   (-1)^{\card{I}}.
$$
Hence
$$
s(\m) = \sum_{j=0}^{h} \binom{j-1}{n-1}  \sum_{I \subseteq [n]:m_{I}=h-j}   (-1)^{\card{I}} 
$$
and, since $\binom{j-1}{n-1} = 0$ whenever $j<n$, we  obtain
$$
s(\m) = \sum_{j=n}^{h} \binom{j-1}{n-1}  \sum_{I \subseteq [n] :m_{I}=h-j}   (-1)^{\card{I}}.
$$
Permuting the order of the two sums  yields
$$
s(\m) =  \sum_{I \subseteq [n] :m_{I} \leq h-n}  \binom{h-m_{I}-1}{n-1} \  (-1)^{\card{I}}  .
$$
\epr

Table~\ref{table.examples.size} illustrates how $s(\m)$ varies as a function of $\m$. 
\begin{table}[hbp]
    \centering
    \begin{tabular}{cr}
        \hline
	$m_{1}, \ldots, m_{n}$ & size  \\
        \hline
	$5, 5$ & 5  \\ 
	$5, 5, 5$ & 19  \\ 
	$5, 5,10$ & 25  \\ 
	$5, 5,100$ & 25  \\ 
	$5, 5,10,10$ & 210   \\
	$5, 5,10,20$ & 250   \\
	$5, 5,10,100$ & 250   \\
	$10,10,10,10$ & 670   \\
	$10,10,10,20$ & 1000   \\
	$10,10,10,100$ & 1000   \\
	$100,100,100,100$ & 666700   \\
	$10, \ldots,10$ (10 times) & 432457640   \\
	$100, \ldots,100$ (10 times) & 430438025018583040   \\
        \hline
    \end{tabular}
    \caption{The size of the largest antichains with $m_{1}, \ldots, m_{n}$ levels on $n$ attributes.}
    \label{table.examples.size}
\end{table}
In particular, we see that increasing one of the $m_{i}$'s way above the others has a limited impact. When all components of $\m$ are identical, it is easy to show that \eqref{eq.Sperner.G} coincides with Sander's expression \eqref{eq:Sander}. In that case, Sander's expression is computationally more efficient than ours.

This section about heterogeneous products does not contain any asymptotic result because it does not seem relevant to let one of the parameters, say $m_{5}$, go to infinity while keeping the other parameters constant.

\section{Homogeneous product}
\label{sec.homo}

Let $h = \floor{n(m+1)/2}$.
Our  result about homogeneous products is the following.

\begin{theo}
\label{theo:Sander.n}
For all $n \geq 2$, if $n(m+1)$ is  even, then $S(m,n)$ is equal to
\begin{multline}
\label{Sander.n.even}
 m^{n-1} 
  - 2  \sum_{r=n-1}^{h-m-1} \ \sum_{i=0}^{\floor{\frac{r-n+1}{m-1}}} (-1)^{i} \binom{n-1}{i} \binom{r-im-1}{r-im-n+1} .  
\end{multline}
Otherwise,  $S(m,n)$ is equal to
\begin{multline}
\label{Sander.n.odd}
 m^{n-1} 
  - 2  \sum_{r=n-1}^{h-m-1} \ \sum_{i=0}^{\floor{\frac{r-n+1}{m-1}}} (-1)^{i} \binom{n-1}{i} \binom{r-im-1}{r-im-n+1} \\
  -  \sum_{i=0}^{\floor{\frac{h-m-n+1}{m-1}}} (-1)^{i} \binom{n-1}{i} \binom{h-im-1}{h-im-n+1}.
\end{multline}
\end{theo}

\bpr
We only prove \eqref{Sander.n.even}. The proof of \eqref{Sander.n.odd} is similar.
We have seen in the proof of Theorem~\ref{theo:SanderGeneralized} that an antichain of maximum size in $\Pi_{i \in [n]} [m_{i}]$ is the set 
$A=\{ x \in X : \sum_{i \in [n]} x_{i} = h \}$ where $h= \floor{\frac{n+\sum_{i \in [n]} m_{i}}{2}}$. Hence, if $n(m+1)$ is even, then an antichain of maximum size in $[m]^{n}$ is the set 
$A=\{ x \in [m]^{n} : \sum_{i \in [n]} x_{i} = h \}$, with $h = n(m+1)/2$.
Since $1 \leq x_{n} \leq m$, 
if  we project the set $A$ on $[m]^{n-1}$ by dropping the last coordinate $x_{n}$, we obtain the set $A' = \{ y \in [m]^{n-1} : h - m \leq \sum_{i \in [n-1]} y_{i} \leq  h-1 \}$. 
Since no $x,y \in A$ are comparable, we know that no distinct $x,y \in A$ project on the same element in $[m]^{n-1}$. Hence $\card{A'} = \card{ A}$ and $S(m,n)$ is equal to
\begin{align*}
\card{ A' } = & \ m^{n-1} - \bcar{ \{ y \in [m]^{n-1} : \sum_{i \in [n-1]} y_{i} \leq  h - m -1 \} } \nonumber \\
  & - \bcar{ \{ y \in [m]^{n-1} :  \sum_{i \in [n-1]} y_{i} \geq  h \} } \nonumber \\
	= & \ m^{n-1} - 2 \bcar{ \{ y \in [m]^{n-1} : \sum_{i \in [n-1]} y_{i} \leq  h - m -1 \} }, 
\end{align*}
where the last equality holds because
$$n(m+1)/2 - m - 1 - \min_{y \in [m]^{n-1}} \sum_{i \in [n-1]} y_{i}  = \max_{y \in [m]^{n-1}} \sum_{i \in [n-1]} y_{i} - h.$$
Let us rewrite $\{ y \in [m]^{n-1} : \sum_{i \in [n-1]} y_{i} \leq  h - m -1 \}$ as the union of several sets:
\begin{equation*}
   \bigcup_{r=n-1}^{h-m-1} A_{r} \text{ where } A_{r} = \{ y \in  [m]^{n-1} :  \sum_{i \in [n-1]} y_{i} = r \}  .\label{eq.union} 
\end{equation*}
Clearly, for any $r \neq s$, $A_{r} \cap A_{s} = \emptyset$ and 
\begin{equation}
S(m,n) = m^{n-1} -2 \sum_{r=n-1}^{h-m-1} \card{A_{r} }.
\label{eq.Smn}
\end{equation}
Let 
\begin{align*}
B_{r} & = \{ y \in  \NN_{+}^{n-1} :  \sum_{i \in [n-1]} y_{i} = r \}, \\
C^{l}_{r} & =  \{y \in \NN_{+}^{n-1}: \sum_{i \in [n-1]} y_{i}= r \text{ and } y_{l} > m \} \text{ for } l \in [n-1]
\end{align*}
and $C^{l*}_{r} =  B_{r} \setminus C^{l}_{r}$.
Then $A_{r} = \bigcap_{l \in [n-1]} C^{l*}_{r}$ and, thanks to the inclusion-exclusion principle, 
\begin{align*}
\card{A_{r}} &= \card{B_{r}} -  \card{\bigcup_{l \in [n-1]} C^{l}_{r} }  \\
 & =   \card{B_{r}} - \sum_{ \emptyset \neq J \subseteq [n-1]} (-1)^{\card{ J} -1} \card{ \bigcap_{l \in J}  C^{l}_{r}} \\
 & =   \card{B_{r}} - \sum_{   [ i ]: 1 \leq i <n} (-1)^{i-1} \binom{n-1}{i} \card{ \bigcap_{l \in [ i ]}  C^{l}_{r}},
\end{align*}
where the last equality holds because all dimensions play the same role. The set $B_{r}$ is a regular $(n-2)$-dimensional simplex. Its cardinality is equal to the $(r-n+1)$-th simplicial polytope number in $n-2$ dimensions \citep{Kim2002}, that is 
\begin{equation}
\label{eq.Br}
 \card{B_{r}}  =  \binom{r-1}{r-n+1}.
 \end{equation}
The set $\bigcap_{l \in [i]}  C^{l}_{r}$ is the set of all elements of $\NN_{+}^{n-1}$ such that at least $i$ components are strictly larger than $m$. If $r < i (m-1)+n-1$, then $\bigcap_{l \in [ i ]}  C^{l}_{r}$ is empty because it is not possible to have at least  $i$ components strictly larger than $m$. Hence 
\begin{equation}
\label{eq.incl.excl}
\card{A_{r}} = \card{B_{r}} - \sum_{   [ i ]: 1 \leq i \leq j} (-1)^{i-1} \binom{n-1}{i} \card{ \bigcap_{l \in [ i ]}  C^{l}_{r}},
\end{equation}
where $j$ is the largest integer such that $r \geq j (m-1)+n-1$.
If $r \geq i (m-1)+n-1$, then 
\begin{equation}
\label{eq.Brim}
\card{ \bigcap_{l \in [ i ]}  C^{l}_{r} } = \card{ B_{r- i m}} =  \binom{r- i m-1}{r- i m-n+1}.
\end{equation}
Combining \eqref{eq.Smn}, \eqref{eq.Br},  \eqref{eq.incl.excl} and \eqref{eq.Brim} concludes the proof.
\epr

\medskip
Expressions~\eqref{Sander.n.even} and~\eqref{Sander.n.odd}  are less elegant than Sander's expression \eqref{eq:Sander}. They are also computationally less efficient. Indeed the first summation in \eqref{Sander.n.even} has approximately $nm/2$ terms while the only summation in \eqref{eq:Sander} has approximately $n/2$ terms. Expressions~\eqref{Sander.n.even} and~\eqref{Sander.n.odd} are nevertheless interesting because they allow us to derive an asymptotic result for $S(m,n)$  when $n$ is fixed and $m \rightarrow \infty$ (see Theorem~\ref{theo:Sander.n.as}). This was not possible with \eqref{eq:Sander}.

When $n < 5$, expressions \eqref{Sander.n.even} and~\eqref{Sander.n.odd} reduce to particularly simple expressions.
\begin{coro}
\label{cor:n<5}
\begin{align*}
S(m,2) &= m; \\
S(m,3) &= \frac{3m^{2}}{4} \text{ if } m \text{ is even and } \frac{3m^{2}+1}{4} \text{ if } m \text{ is odd}; \\
S(m,4) &= \frac{2m^{3}+m}{3}.
\end{align*}
\end{coro}
For $n =2, 3$ or $4$, the asympotic behaviour of $S(m,n)$  is easy to derive from this corollary, while the general case is covered by our next result.

\medskip

\begin{theo}
\label{theo:Sander.n.as}
For all $n \geq 2$, when $m \rightarrow \infty$,   $S(m,n)$ is equal to $m^{n-1} g(n) + \bigO(m^{n-2})$ where $g(n)$ is equal to 
\begin{equation}
\label{Sander.n.even.as}
  1 - 2 \sum_{j=0}^{\frac{n-4}{2}} \ \sum_{i=0}^{j}  (-1)^{i}  \ \frac{(1+j-i)^{n-1}-(j-i)^{n-1} }{i! \ (n-1-i)!},
\end{equation}
when $n$ is even, 
or to 
\begin{multline}
\label{Sander.n.odd.as}
  1 - 2 \sum_{j=0}^{\frac{n-5}{2}} \ \sum_{i=0}^{j}  (-1)^{i}  \ \frac{(1+j-i)^{n-1}-(j-i)^{n-1} }{i! \ (n-1-i)!} \\
  - 2 \sum_{i=0}^{\frac{n-3}{2}} \  (-1)^{i}  \ \frac{(\frac{n}{2}-1-i)^{n-1}-(\frac{n-3}{2}-i)^{n-1} }{i! \ (n-1-i)!},
\end{multline}
when $n$ is odd.
\end{theo}

\bpr
We only prove \eqref{Sander.n.even.as}.
Expression~\eqref{Sander.n.even} for $S(m,n)$ can also be written as
\begin{multline}
\label{Sander.n.even.alt}
 m^{n-1}   - 2 \left( 
\sum_{j=0}^{\frac{n-4}{2}} \ \sum_{r=n-1+j(m-1)}^{m+n-3+j(m-1)} \ \sum_{i=0}^{j} (-1)^{i} \binom{n-1}{i} \binom{r-im-1}{r-im-n+1}
\right.  \\
  + \left. \sum_{r=n-1+\frac{n-2}{2}(m-1)}^{h - m -1} \ \sum_{i=0}^{\frac{n-2}{2}} (-1)^{i} \binom{n-1}{i} \binom{r-im-1}{r-im-n+1} \right) .  
\end{multline}
For $n$ fixed, $g(n)$ is the limit for $n \rightarrow \infty$ of \eqref{Sander.n.even.alt} divided by $m^{n-1}$, that is
\begin{multline*}
 1  - \lim_{m \rightarrow \infty} \frac{2}{(n-2)! \, m^{n-1}} \left( 
\sum_{j=0}^{\frac{n-4}{2}} \ \sum_{i=0}^{j} (-1)^{i} \binom{n-1}{i} \ \sum_{r=n-1+j(m-1)}^{m+n-3+j(m-1)}   \frac{(r-im-1)!}{(r-im-n+1)!}
\right.  \\
  + \left.  \sum_{i=0}^{\frac{n-2}{2}} (-1)^{i} \binom{n-1}{i}  \sum_{r=n-1+\frac{n-2}{2}(m-1)}^{h-m-1} \frac{(r-im-1)!}{(r-im-n+1)!}\right) .
\end{multline*}
Since $\lim_{a \rightarrow \infty} \frac{a!/(a-b)!}{a^{b}} = 1$, this is also equal to 
\begin{multline*}
 1  - \lim_{m \rightarrow \infty} \frac{2}{(n-2)! \, m^{n-1}} \left( 
\sum_{j=0}^{\frac{n-4}{2}} \ \sum_{i=0}^{j} (-1)^{i} \binom{n-1}{i} \ \sum_{r=n-1+j(m-1)}^{m+n-3+j(m-1)}   (r-im)^{n-2}
\right.  \\
  + \left.  \sum_{i=0}^{\frac{n-2}{2}} (-1)^{i} \binom{n-1}{i}  \sum_{r=n-1+\frac{n-2}{2}(m-1)}^{h-m-1}  (r-im)^{n-2} \right) 
\end{multline*}
or, since $\sum_{d=a}^{b} d^{c} = (c+1)^{-1}[d^{c+1}]_{a}^{b} + \mathcal{O}(d^{c})$ as $d \rightarrow \infty$,
\begin{multline*}
 1 - \lim_{m \rightarrow \infty} \frac{2}{(n-1)! \, m^{n-1}} \left( 
\sum_{j=0}^{\frac{n-4}{2}} \ \sum_{i=0}^{j} (-1)^{i} \binom{n-1}{i}   \Big[(r-im)^{n-1}\Big]_{r=n-1+j(m-1)}^{m+n-3+j(m-1)}
\right.  \\
  + \left.  \sum_{i=0}^{\frac{n-4}{2}} (-1)^{i} \binom{n-1}{i}   \Big[(r-im)^{n-1}\Big]_{r=n-1+\frac{(n-2)}{2}(m-1)}^{h-m-1} \right) .
\end{multline*}
We then substitute $r$ with the summation bounds, we simplify and we take the limit, keeping only the coefficients of the highest power of $m$ (i.e.\ $m^{n-1}$) and we obtain \eqref{Sander.n.even.as}. The second highest power of $m$ is $m^{n-2}$ and this completes the proof.
\epr

Expressions \eqref{Sander.n.even.as} and \eqref{Sander.n.odd.as} are easy to compute for $n$ between 2 and approximately 100. Beyond 100, using Stirling's approximation for  factorials helps. For very large $n$, the computation time becomes prohibitive because of the large number of terms in the summations.

\begin{figure}
\begin{center}
\scalebox{0.5}{\includegraphics{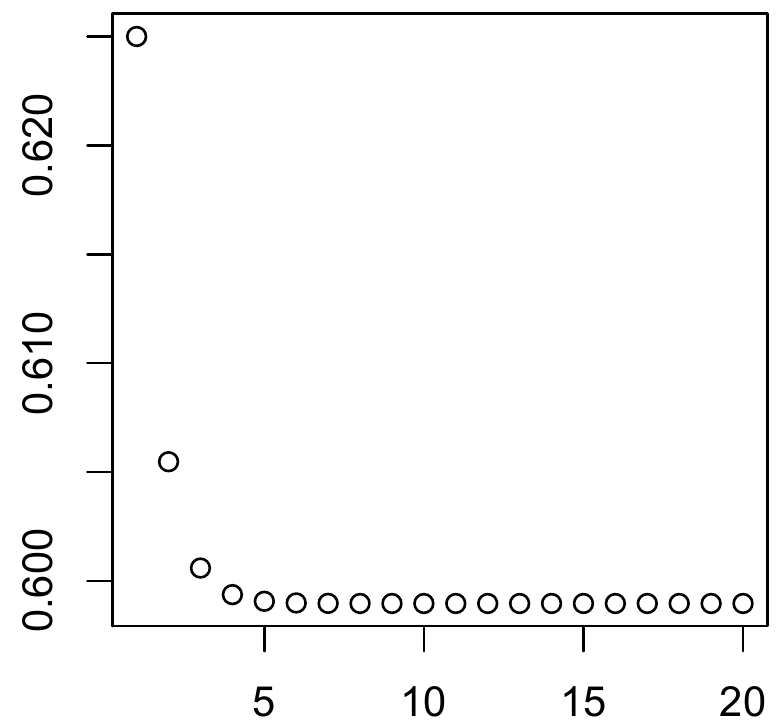}} \scalebox{0.5}{\includegraphics{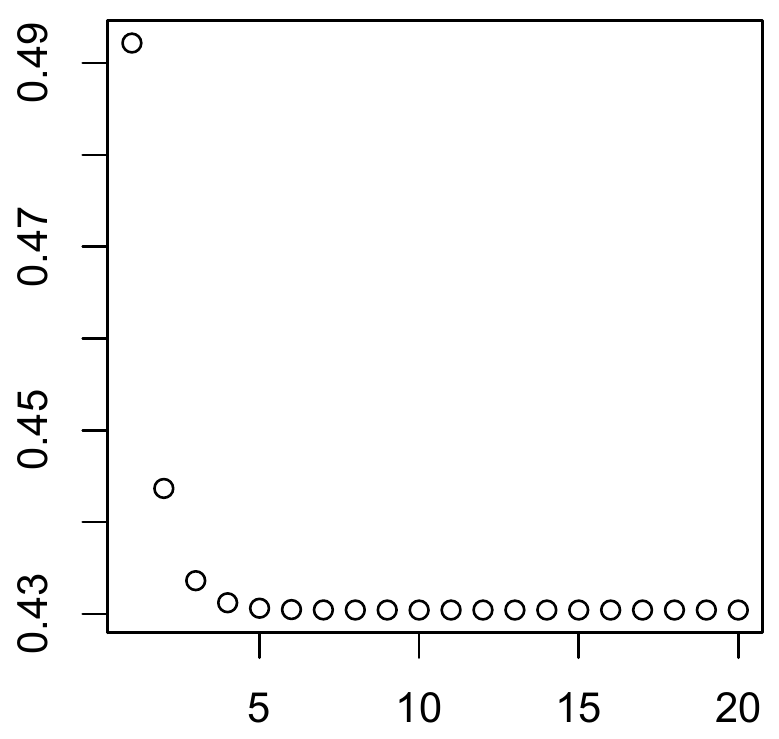}}
\caption{Horizontal axis: $\log_{2}m$. Vertical axis: ratio between $S(m,n)$ and $m^{n-1}$ for various values of  $m$. Left: $n=5$; right: $n=10$.}
\label{fig.convergence}
\end{center}
\end{figure}

Numerical estimations of  $g(n)$ are given in Table~\ref{table.ratio} and
Figure~\ref{fig.convergence} illustrates how quickly $S(m,n)/m^{n-1}$ converges to $g(n)$ for $n=5$ and $10$.
\begin{table}[hbp]
    \centering
    \begin{tabular}{llll}
        \hline
        $n$ & $g(n)$ & $n$ & $g(n)$  \\
        \hline
        2 & 1 & 8 & 0.47936507936507944 \\
        3 & 3/4 & 9 & 0.45292096819196426 \\
        4 & 2/3  & 10 & 0.43041776895943573 \\
        5 & 115/192 & 20 & 0.30669310173797590  \\
        6 & 11/20 & 30 & 0.25104851499027436  \\
        7 & 0.5110243055555556  & 100 & 0.12337408033008801 \\
          \hline
    \end{tabular}
    \caption{Some values of $g(n)$. Fractions are exact.}
    \label{table.ratio}
\end{table}

Notice that, when $n < 5$, Corollary~\ref{cor:n<5}  provides an asymptotic expression for $S(m,n)$ that is tighter than that resulting from Theorem~\ref{theo:Sander.n.as}.
 
 
 \section*{References}

 \bibliography{mybib}

\begin{thebibliography}{15}
\expandafter\ifx\csname natexlab\endcsname\relax\def\natexlab#1{#1}\fi
\expandafter\ifx\csname url\endcsname\relax
  \def\url#1{\texttt{#1}}\fi
\expandafter\ifx\csname urlprefix\endcsname\relax\def\urlprefix{URL }\fi

\bibitem[{Caspard et~al.(2012)Caspard, Leclerc, and
  Monjardet}]{CaspardLeclercMonjardet2012}
Caspard, N., Leclerc, B., Monjardet, B., 2012. Finite Ordered Sets: Concepts,
  Results and Uses. Encyclopedia of Mathematics and its Applications. Cambridge
  University Press.

\bibitem[{de~{M}oivre(1756)}]{Moivre1756}
de~{M}oivre, A., 1756. The doctrine of chances. Third ed., reprinted by
  Chelsea, New York, 1967.

\bibitem[{{Ersek Uyan{\i}k} et~al.(2017){Ersek Uyan{\i}k}, Sobrie, Mousseau,
  and Pirlot}]{ErsekUyanikSMP17}
{Ersek Uyan{\i}k}, E., Sobrie, O., Mousseau, V., Pirlot, M., 2017. Enumerating
  and categorizing positive {Boolean} functions separable by a $k$-additive
  capacity. Discrete Applied Mathematics 229, 17--30.

\bibitem[{Fern\'andez et~al.(2017)Fern\'andez, Figueira, Navarro, and
  Roy}]{FernandezFigueiraNavarroRoy2017}
Fern\'andez, E., Figueira, J.~R., Navarro, J., Roy, B., 2017. {ELECTRE}
  {TRI-nB}: A new multiple criteria ordinal classification method. European
  Journal of Operational Research 263~(1), 214--224.

\bibitem[{Grabisch(2016)}]{Grabisch2016}
Grabisch, M., 2016. Remarkable polyhedra related to set functions, games and
  capacities. TOP 24~(2), 301--326.

\bibitem[{Hsiao and Raghavan(1993)}]{HsiaoRaghavan1993}
Hsiao, C.-R., Raghavan, T. E.~S., 1993. Shapley value for multichoice
  cooperative games, {I}. Games and Economic Behavior 5~(2), 240 -- 256.

\bibitem[{Kim(2002)}]{Kim2002}
Kim, H.~K., 2002. On regular polytope numbers. Proceedings of the American
  Mathematical Society 131~(1), 65--75.

\bibitem[{{MacMahon}(1893)}]{MacMahon1893}
{MacMahon}, P.~A., 1893. Memoir on the theory of the compositions of numbers.
  Philosophical Transactions of the Royal Society A 184, 835--901.

\bibitem[{Mattner and Roos(2008)}]{MattnerRoos2008}
Mattner, L., Roos, B., 2008. Maximal probabilities of convolution powers of
  discrete uniform distributions. Statistics and Probability Letters 78,
  2992--2996.

\bibitem[{Motek(1986)}]{Motek1986}
Motek, J., 1986. Problem 86-8. The Mathematical Intelligencer 8.

\bibitem[{OEIS(2019)}]{OEISB}
OEIS, 2019. The On-line Encyclopaedia of Integer Sequences, {\emph{Sloane, N.
  J. A. (Ed.)}}.
\newline\urlprefix\url{https://oeis.org}

\bibitem[{Proctor et~al.(1980)Proctor, Saks, and
  Sturtevant}]{ProctorSaksSturtevant1980}
Proctor, R.~A., Saks, M.~E., Sturtevant, D.~G., 1980. Product partial orders
  with the {S}perner property. Discrete Mathematics 30, 173--180.

\bibitem[{Sander(1993)}]{Sander1993}
Sander, J.~W., 1993. On maximal antihierarchic sets of integers. Discrete
  Mathematics 113, 179--189.

\bibitem[{Sperner(1928)}]{Sperner1928}
Sperner, E., 1928. Ein {S}atz {\"u}ber {U}ntermengen einer endlichen {M}enge.
  Mathematische Zeitschrift 27~(1), 544--548.

\bibitem[{Tsai(2018)}]{Tsai2018}
Tsai, S.-F., Dec 2018. A simple upper bound on the number of antichains in
  $[t]^n$. Order.
\newline\urlprefix\url{https://doi.org/10.1007/s11083-018-9480-5}

\end{thebibliography}

\end{document}